\documentclass[10pt,twoside]{article}
\usepackage{graphicx}
\usepackage{amsmath}
\usepackage{amscd,amsmath,amssymb,amsfonts}
\usepackage{amssymb}
\usepackage{Latex-document}

\title{\bf  Toeplitz Determinants, Random Growth \vskip -2mm
and Determinantal Processes\vskip 6mm}

\author{K. Johansson\vspace*{-0.5cm}\thanks{Department of Mathematics,
Royal Institute of Technology, SE-100 44 Stockholm, Sweden.
E-mail: kurtj@math.kth.se}}
\date{}

\newcommand{\Ai}{\text{Ai\,}}

\begin{document}
\maketitle

\thispagestyle{first} \setcounter{page}{53}

\begin{abstract}

\vskip 3mm

We summarize some of the recent developments which link certain
problems in combinatorial theory related to random growth to random
matrix theory.

\vskip 4.5mm


\end{abstract}

\vskip 12mm

\section{Introduction} \label{section 1}\setzero

\vskip -5mm \hspace{5mm}

Let $\sigma$ be a permutation from $S_N$. We say that
$\sigma(i_1),\dots,\sigma(i_m)$, $i_1<\dots<i_m$, is an {\it increasing
subsequence} of $\sigma$ if $\sigma(i_1)<\dots<\sigma(i_m)$. The number
$m$ is the length of the subsequence. The length of the longest
increasing subsequence in $\sigma$ is denoted by $\ell_N(\sigma)$. If
we pick $\sigma$ from $S_N$ uniformly at random $\ell_N(\sigma)$
becomes a random variable. {\it Ulam's problem}, \cite{Ul}, is the study
of the asymptotic properties as $N\to\infty$
of this random variable in particular its mean. It turns out that
there is a surprisingly rich mathematical structure around this
problem as we hope will be clear from the presentation below. It has
been known for some time that $\mathbb{E}[\ell_N]\sim 2\sqrt{N}$
as $N\to\infty$, \cite{VeKe}, \cite{LoSh}. We refer to \cite{AlDi} for
some background to the problem. A Poissonized version of the problem
can be obtained by letting $N$ be an independent Poisson random
variable with mean $\alpha$. This gives a random variable $L(\alpha)$
with distribution
\begin{equation}\label{1.1}
\mathbb{P}[L(\alpha)\le
n]=\sum_{N=0}^\infty\frac{e^{-\alpha}\alpha^N}{N!}
\mathbb{P}[\ell_N\le n].
\end{equation}
Since $\mathbb{P}[\ell_N\le n]$ is a decreasing function of $N$,
\cite{Jo1}, asymptotics of the left hand side of (1.1) can be used to
obtain asymptotics of $\mathbb{P}[\ell_N\le n]$ (de-Poissonization).

The random variable $L(\alpha)$ can be realized geometrically using
Hammersley's picture, \cite{Ha}. Consider a Poisson process with
intensity 1 in the square $[0,\gamma]^2$, $\gamma=\sqrt{\alpha}$. An
{\it up/right path} is a sequence of Poisson points
$(x_1,y_1),\dots, (x_m,y_m)$ in the square such that $x_i<x_{i+1}$ and
$y_i<y_{i+1}$, $i=1,\dots,m-1$. The maximal number of points
in an up/right path has the same distribution as $L(\alpha)$. A
sequence of points realizing this maximum is called a {\it maximal path}. It
is expected from heuristic arguments, see below, that the standard
deviation of $L(\alpha)$ should be of order
$\gamma^{1/3}=\alpha^{1/6}$. The proof that this is true, \cite{BDJ},
and that we can also understand the law of the fluctuations is the
main recent result that
will be discussed below. Also, the deviations of a maximal path from
the diagonal $x=y$ should be of order $\gamma^{2/3}$. This last
statement is proved in \cite{Jo3}.

A generalization of the random variable $L(\alpha)$ can be defined in
the following way. Let $w(i,j)$, $(i,j)\in\mathbb{Z}_+^2$, be
independent geometric random variables with parameter $q$. An {\it up/right
path} $\pi$ from $(1,1)$ to $(M,N)$ is a sequence
$(1,1)=(i_1,j_1),(i_2,j_2),\dots ,(i_m,j_m)=(M,N)$, $m=M+N-1$, such
that either $i_{r+1}-i_r=1$ and $j_{r+1}=j_r$, or
$i_{r+1}=i_r$ and $j_{r+1}-j_r=1$. Set
\begin{equation}\label{1.2}
G(M,N)=\max_{\pi}\sum_{(i,j)\in\pi} w(i,j),
\end{equation}
where the maximum is taken over all up/right paths $\pi$ from $(1,1)$
to $(M,N)$. Alternatively, we can define $G(M,N)$ recursively by
\begin{equation}\label{1.3}
G(M,N)=\max(G(M-1,N),G(M,N-1))+w(M,N).
\end{equation}
Some thought shows that if we let $q=\alpha/N^2$ then $G(N,N)$
converges in distribution to $L(\alpha)$ as $N\to\infty$,
\cite{Jo4}, so we can view
$G(N,N)$ as a generalization of $L(\alpha)$. We can think of
(\ref{1.2}) as a directed last-passage site percolation problem. Since
all paths $\pi$ have the same length, if $w(i,j)$ were a bounded
random variable we could relate (\ref{1.2}) to the corresponding
first-passage site percolation problem, with a min instead of a max. in
(\ref{1.2}). The random variable $G(M,N)$ connects with many different
problems, a corner growth model, zero-temperature directed polymers,
totally asymmetric simple exclusion processes and domino tilings of
the Aztec diamond, see \cite{Jo2}, \cite{Jo4} and references
therein. It is also related to another growth model, the {\it (discrete)
polynuclear growth (PNG) model}, \cite{KrSp}, \cite{PrSp} defined as
follows. Let $h(x,t)\in\mathbb{N}$ denote the height above
$x\in\mathbb{Z}$ at time $t\in\mathbb{N}$. The growth model is defined
by the recursion
\begin{equation}\label{1.4}
h(x,t+1)=\max (h(x-1,t), h(x,t),h(x+1,t)) +a(x,t),
\end{equation}
where $a(x,t)$, $(x,t)\in\mathbb{Z}\times\mathbb{N}$, are independent
random variables. If we assume that $a(x,t)=0$ whenever $x-t$ is even, and
that the distribution of $a(x,t)$ is geometric with parameter $q$,
then setting $w(i,j)=a(i-j,i+j-1)$, we obtain
$G(i,j)=h(i-j,i+j-1)$. The growth model (\ref{1.4}) has some relation
to the so called Kardar-Parisi-Zhang (KPZ) equation, \cite{KrSp}, and
is expected to fall within the so called KPZ-universality class. The
exponents 1/3 and 2/3 discussed above are the conjectured exponents
for $1+1$-dimensional growth models in this class.

\section{Orthogonal polynomial ensembles} \label{section 2}

\vskip -5mm \hspace{5mm}

Consider a probability (density) on $\Omega^N$,
$\Omega=\mathbb{R}$, $\mathbb{Z}$, $\mathbb{N}$ or
$\{0,1,\dots,M\}$ of the form
\begin{equation}\label{2.1}
u_N(x)=\frac 1{Z_N}\Delta_N(x)^2\prod_{j=1}^M w(x_j),
\end{equation}
where $\Delta_N(x)=\prod_{1\le i<j\le N}(x_i-x_j)$ is the Vandermonde
determinant, $w(x)$ is some non-negative weight function on $\Omega$
and $Z_N$ is a normalization constant. We call such a probability an
{\it orthogonal polynomial ensemble}. We can think of this as a finite
point process on $\Omega$. Let $d\mu$ be Lebesgue or counting measure
on $\Omega$ and let $p_n(x)=\kappa_nx^n+\dots$ be the normalized
orthogonal polynomials with respect to the measure $w(x)d\mu(x)$ on
$\Omega$. The correlation functions $\rho_{m,N}(x_1,\dots,x_m)$ of the
point process are given by determinants, we have a so called
{\it determinantal point process}, \cite{So}. In fact,
\begin{align}\label{2.2}
\rho_{m,N}(x_1,\dots,x_m)&=\frac {N!}{(N-m)!} \int_{\Omega^{N-m}}
u_N(x) d\mu(x_{m+1})\dots d\mu(x_N)\\
&=\det(K_N(x_i,x_j))_{1\le i,j\le m},
\notag
\end{align}
where the kernel $K_N$ is given by
\begin{equation}\label{2.3}
K_N(x,y)=\frac{\kappa_{N-1}}{\kappa_N}
\frac{p_N(x)p_{N-1}(y)-p_{N-1}(x)p_N(y)}{x-y} (w(x)w(y))^{1/2}.
\end{equation}
A computation shows that for bounded $f:\Omega\to\mathbb{C}$,
\begin{align}\label{2.4}
\mathbb{E}[\prod_{j=1}^N(1+f(x_j))]&=\sum_{k=0}^N\frac 1{k!}
\int_{\Omega^k} \prod_{j=1}^k f(x_j)
\det(K_N(x_i,x_j))_{1\le i,j\le k}d^k\mu(x)\\
&=\det (I+fK_N)_{L^2(\Omega,d\mu)},
\notag
\end{align}
where the last determinant is the Fredholm determinant of the integral
operator on $L^2(\Omega, d\mu)$ with kernel $f(x)K_N(x,y)$. In
particular, we can compute hole or gap probabilities, e.g. the
probability of having no particle in an interval $I\subseteq\Omega$
by taking $f=-\chi_I$, minus the characteristic function of the
interval $I$. If $x_{\text{max}}=\max x_j$ denotes the position of the
rightmost particle it follows that
\begin{equation}\label{2.5}
\mathbb{P}[x_{\text{max}}\le a]=\det(I-K_N)_{L^2((a,\infty),d\mu)}.
\end{equation}

As $N\to\infty$ we can obtain limiting determinantal processes on
$\mathbb{R}$ or $\mathbb{Z}$ with kernel $K$, i.e. the probability
(density) of finding particles at $x_1,\dots, x_m$ is given by
$\det(K(x_i,x_m))_{1\le i,j\le m}$. We will be interested in the limit
process around the rightmost particle. This is typically given by the
{\it Airy kernel}, $x,y\in\mathbb{R}$,
\begin{equation}\label{2.6}
A(x,y)=\frac{\Ai (x)\Ai'(y)-\Ai'(x)\Ai(y)}{x-y}.
\end{equation}
This limit, the {\it Airy point process} has a rightmost particle
almost surely and its position has the distribution function
\begin{equation}\label{2.7}
F_2(\xi)=\det(I-A)_{L^2(\xi,\infty)},
\end{equation}
known as the {\it Tracy-Widom distribution}, \cite{TW1}. Hole
probabilities as functions of the endpoints of the intervals
satisfy systems of differential equations, \cite{TW2},
\cite{ASvM}. For example, we have
\begin{equation}\label{2.8}
F_2(\xi)=\exp(-\int_{\xi}^\infty (x-\xi)u(x)^2dx),
\end{equation}
where $u$ solves the Painlev\'e-II equation $u''=xu+2u^3$ with
boundary condition $u(x)\sim Ai(x)$ as $x\to\infty$.

An example of a measure of the form (\ref{2.1}) comes from the {\it
Gaussian Unitary Ensemble} ({\it GUE}) of random matrices. The GUE
is a Gaussian measure on the space $\mathcal{H}_N\cong
\mathbb{R}^{N^2}$ of all $N\times N$ Hermitian matrices. It is defined
by $d\mu_{\text{GUE},N}(M)=Z_N^{-1}\exp(-\text{tr\,} M^2)dM$, where
$dM$ is the Lebesgue measure on $\mathcal{H}_N$ and $Z_N$ is a
normalization constant. The corresponding eigenvalue measure has the
form (\ref{2.1}) with $w(x)=\exp(-x^2)$ and $\Omega=\mathbb{R}$,
\cite{Me}. Hence the $p_n$:s are multiples of the ordinary Hermite
polynomials. The largest eigenvalue $x_{\max}$ will lie around
$\sqrt{2N}$. This is related to the fact that the largest zero of
$p_N$ lies around $\sqrt{2N}$. The local asymptotics of
$p_N(x)\exp(-x^2/2)$ around this point,
$x=\sqrt{2N}+\xi/N^{1/6}\sqrt{2}$, is given by the Airy function,
$\Ai(\xi)$. This asymptotics, some estimates, (\ref{2.5}) and
(\ref{2.7}) give the following result,
\begin{equation}\label{2.9}
\mathbb{P}_{\text{GUE},N}[\frac {\sqrt{2N}x_{\max} -2N}{N^{1/3}}\le
\xi]\to F_2(\xi)
\end{equation}
as $N\to\infty$.

\section{Some theorems} \label{section 3}

\vskip -5mm \hspace{5mm}

The previous section may seem unrelated to the first but as the next
theorems will show the problem of understanding the distribution of
$L(\alpha)$ and $G(M,N)$ fits nicely into the machinery of sect. 2.

{\bf Theorem 3.1.} \cite{Jo2}. \it  Take $\Omega=\mathbb{N}$, $M\ge N$
and $w(x)=\binom{M-N+x}{x}$ in (\ref{2.1}). Then $G(M,N)$ is
distributed exactly as $x_{\max}$. \rm

The corresponding orthogonal polynomials are the Meixner polynomials,
a classical family of discrete orthogonal polynomials, and we refer to
the measure obtained as the {\it Meixner ensemble}. It is an example of a
{\it discrete orthogonal polynomial ensemble}, \cite{Jo4}. By computing the
appropriate Airy asymptotics of the Meixner polynomials we can use
(\ref{2.5}) to prove the next theorem.

{\bf Theorem 3.2.} \cite{Jo2}. \it Let $\gamma\ge 1$ be fixed and set
$\omega(\gamma,q)=(1-q)^{-1}(1+\sqrt{q\gamma})^2-1$ and $\sigma(\gamma,q)=
(1-q)^{-1}(q/\gamma)^{1/6}(\sqrt{\gamma}+\sqrt{q})^{2/3}
(1+\sqrt{q\gamma})^{2/3}$. Then,
\begin{equation}\label{3.1}
\lim_{N\to\infty}\mathbb{P}[\frac{G([\gamma N],N)-\omega(\gamma,q)N}{\sigma(\gamma,q)
  N^{1/3}} \le\xi]=F_2(\xi).
\end{equation}
\rm

Thus $G([\gamma N],N)$ fluctuates like the largest eigenvalue of a GUE
matrix.

As discussed above by setting $q=\alpha/N^2$ we can obtain $L(\alpha)$
as a limit of $G(N,N)$ as $N\to\infty$. By taking this limit in
theorem 3.1, and using the fact that the measure has determinantal
correlation functions, we see that $L(\alpha)$ behaves like the
rightmost particle in a determinantal point process on $\mathbb{Z}$
given by the {\it discrete Bessel kernel}, \cite{Jo4}, \cite{BOO},
\begin{equation}\label{3.2}
B^\alpha(x,y)=\sqrt{\alpha}\frac{J_x(2\sqrt{\alpha})J_{y+1}(2\sqrt{\alpha})
-J_{x+1}(2\sqrt{\alpha})J_y(2\sqrt{\alpha})}{x-y},
\end{equation}
$x,y\in\mathbb{Z}$. This gives
\begin{equation}\label{3.3}
\mathbb{P}[L(\alpha)\le n]=\det(I-B^\alpha)_{\ell^2(\{n,n+1,\dots\})}.
\end{equation}
Once we have this formula we see that all we need is
the classical asymptotic formula $\alpha^{1/6}
J_{2\sqrt{\alpha}+\xi\alpha^{1/6}}(2\sqrt{\alpha}
\to\Ai (\xi)$ as $\alpha\to\infty$
uniformly in compact intervals, and some estimates of the
Bessel functions in order to get a limit theorem for $L(\alpha)$:

{\bf Theorem 3.3.} \cite{BDJ}. \it As $\alpha \to\infty$,
\begin{equation}\label{3.4}
\mathbb{P}[\frac{L(\alpha)-2\sqrt{\alpha}}{\alpha^{1/6}}\le\xi ]\to F_2(\xi).
\end{equation}
\rm

Note the similarity with (\ref{2.9}), just replace $N$ by
$\sqrt{\alpha}$. This result was first proved in \cite{BDJ} by another
method, see below. De-poissonizing we get a limit theorem for
$\ell_N(\pi)$, see \cite{BDJ}.

\section{Rewriting Toeplitz determinants} \label{section 4}

\vskip -5mm \hspace{5mm}

The {\it Toeplitz determinant} of order $n$ with generating function
$f\in L^1(\mathbb{T})$ is defined by
\begin{equation}\label{4.1}
D_n(f)=\det(\hat{f}_{i-j})_{1\le i,j\le n},
\end{equation}
where $\hat{f}_k=(2\pi)^{-1}\int_{-\pi}^\pi
f(e^{i\theta})e^{-ik\theta}d\theta$ are the complex Fourier
coefficients of $f$. Consider the generating function
\begin{equation}\label{4.2}
f(z)=\prod_{\ell=1}^M(1+\frac {a_\ell}z)(1+b_\ell z),
\end{equation}
where $a_\ell,b_\ell$ are complex numbers. The elementary symmetric
polynomial $e_m(a)$, $a=(a_1,\dots,a_M)$ is defined by $\prod_{j=1}^M(1+a_j
z)=\sum_{|m|<\infty} e_m(a)z^m$. A straightforward computation shows
that when $f$ is given by (\ref{4.2}) then
\begin{equation}
\hat{f}_{i-j}=\sum_{m=0}^\infty e_{m-j}(a)e_{m-k}(b).
\notag
\end{equation}
Insert this into the definition (\ref{4.1}) and use the Heine
identity,
\begin{align}\label{4.3}
\frac 1{n!}\int_{\Omega^n} &\det(\phi_i(x_j))_{1\le i,j\le n}
\det(\psi_i(x_j))_{1\le i,j\le n}d^n\mu(x)\\
&=\det(\int_\Omega \phi_i(x)\psi_j(x) d\mu(x))_{1\le i,j\le n},
\notag
\end{align}
to see that
\begin{equation}\label{4.4}
D_n(f)=\sum_{m_1>m_2>\dots>m_n\ge 0} \det(e_{m_i-j}(a))_{1\le i,j\le n}
\det(e_{m_i-j}(b))_{1\le i,j\le n}.
\end{equation}
Here we have removed the $n!$ by ordering the variables. These
determinants are again symmetric polynomials, the so called Schur
polynomials. Let $\lambda=(\lambda_1,\lambda_2,\dots)$ be a partition
and let
$\lambda'=(\lambda_1',\lambda_2',\dots)$ be the conjugate partition,
\cite {Sa}. Set $m_i=\lambda_i'+n-i$, $i=1,\dots,n$ and $\lambda_i'=0$
if $i>n$, so that $\lambda'$ has at most $n$ parts, $\ell(\lambda')\le
n$, which means that $\lambda_1\le n$. Then the Schur polynomial
$s_\lambda(a)$ is given by
\begin{equation}\label{4.5}
s_\lambda(a)=\det(e_{\lambda_i'-i+j}(a))_{1\le i,j\le n}=
\det(e_{m_i-j}(a))_{1\le i,j\le n},
\end{equation}
the Jacobi-Trudi identity. Hence,
\begin{equation}\label{4.6}
D_n(f)=\sum_{\lambda\,;\,\lambda_1\le n} s_\lambda(a)s_\lambda(b),
\end{equation}
and we have derived Gessel's formula, \cite{Ge}. If we let
$n\to\infty$ in the right hand side we obtain
$\prod_{i,j=1}^M(1-a_ib_j)^{-1}$ by the Cauchy identity, \cite{Sa}. In
the case when all $a_i,b_j\in[0,1]$, $s_\lambda(a)s_\lambda(b)\ge 0$,
and we can think of
\begin{equation}\label{4.7}
 \prod_{i,j=1}^M(1-a_ib_j) s_\lambda(a)s_\lambda(b)
\end{equation}
as a probability measure on all partitions $\lambda$ with at most $n$
parts, the {\it Schur measure}, \cite{Ok2}. In this formula we can
insert the combinatorial definition of the Schur polynomial, \cite
{Sa},
\begin{equation}\label{4.8}
s_\lambda(a)=\sum_{T\,:\,\text{sh\,}(T)=\lambda} a_1^{m_1(T)}\dots
a_M^{m_M(T)},
\end{equation}
where the sum is over all semi-standard Young tableaux $T$,
\cite{Sa}, with shape $\lambda$, and $m_i(T)$ is the number of $i$:s
in $T$.

A connection with the random variables in section 1 is now provided by
the Robinson-Schensted-Knuth (RSK) correspondence, \cite{Sa}. This
correspondence maps an $M\times M$ integer matrix to a pair of
semi-standard Young tableaux $(T,S)$ with entries from $\{1,2,\dots,
M\}$. If we let the random variables $w(i,j)$ be independent geometric
with parameter $a_ib_j$ then the RSK-correspondence maps the measure
we get on the integer matrix $(w(i,j))_{1\le i,j\le M}$ to the Schur
measure (\ref{4.7}).
Also, the RSK-correspondence is such that $G(M,M)=\lambda_1$, the
length of the first row. If we put $a_j=0$ for $N<j\le M$ and
$a_i=b_i=\sqrt{q}$ for $1\le i\le N$ in the Schur measure and set
$x_j=\lambda_j+N-j$, $1\le j\le N$, we obtain the result in Theorem
3.1.

In the limit $M=N\to\infty$, $q=\alpha/N^2$, in which case $G(N,N)$
converges to $L(\alpha)$, the Schur measure converges to the so called
{\it Plancherel measure} on partitions, \cite{VeKe},
\cite{Jo4}. In the variables
$\lambda_i-i$ this measure is a determinantal point process on
$\mathbb{Z}$ given by the kernel $B^\alpha$, (\ref{3.2}). This result
was obtained independently in \cite{BOO}, which also gives a
descrption in terms of different coordinates. See also \cite{Ok1} for a
more direct geometric relation between GUE and the Plancherel
measure. In this limit the Toeplitz determinant formula (\ref{4.6})
gives
\begin{equation}\label{4.9}
\mathbb{P}[L(\alpha)\le
n]=e^{-\alpha}D_n(e^{2\sqrt{\alpha}\cos\theta}).
\end{equation}
This variant of Gessel's formula was the starting point for the
original proof of Theorem 3.3 in \cite{BDJ}. The right hand side of
(\ref{4.9}) can be expressed in terms of the leading coefficients of
the orthogonal polynomials on $\mathbb{T}$ with respect to the weight
$\exp(2\sqrt{\alpha}\cos\theta)$. These orthogonal polynomials in turn
can be obtained as a solution to a matrix-valued Riemann-Hilbert
problem (RHP), and the asymptotics of this RHP as $\alpha\to\infty$
can be analyzes using the powerful asymptotic techniques developed by
Deift and Zhou, \cite{DKMVZ}. This approach leads to the formula
(\ref{2.8}) for the limiting distribution.

Write $f=\exp(g)$ and insert the definition of the Fourier
coefficients into the definition (\ref{4.1}). By the Heine identity we
obtain an integral formula for the Toeplitz determinant,
\begin{align}\label{4.10}
D_n(f)&=\frac 1{(2\pi)^nn!}\int_{[-\pi,\pi]^n} \prod_{1\le\mu<\nu\le n}
|e^{i\theta_\mu}- e^{i\theta_\nu}|^2 \prod_{\mu=1}^n
e^{g(e^{i\theta_\mu})} d^n\theta\\
&=\int_{U(n)} e^{\text{tr\,}g(U)} dU.
\notag
\end{align}
In the last integral $dU$ denotes normalized Haar measure on the
unitary group $U(n)$ and the identity is the Weyl integration
formula. The limit of (\ref{4.9}) as $\alpha\to\infty$ is then a so
called {\it double scaling limit} in a unitary matrix model,
\cite{PeSh}. The formula (\ref{4.9}) can also be obtained by
considering the integral over the unitary group, see \cite{Ra}.

Another way to obtain the Schur measure is via families of
non-intersecting paths which result from a multi-layer PNG model,
\cite{Jo5}. The determinants in the measure then come from the
Karlin-McGregor theorem or the Lindstr\"om-Gessel-Viennot method.

\section{A curiosity} \label{section 5}

\vskip -5mm \hspace{5mm}

Non-intersecting paths can also be used to describe certain tilings,
e.g. domino tilings and tilings of a hexagon by rhombi. By looking at
intersections with appropriate lines one can obtain discrete
orthogonal polynomial ensembles. In the case of tilings of a hexagon
by rhombi, which correspond to boxed planar partitions, \cite{CLP}, the
Hahn ensemble, i.e. (\ref{2.1}) with $\Omega=\{0,\dots,M\}$ and a
weight giving the Hahn polynomials, is obtained, \cite{Jo5}. The
computation leading to this result also gives a proof of the classical
MacMahon formula, \cite{St}, for the number of boxed planar partitions
in an $abc$ cube, i.e. the number of rhombus tilings of an
$abc$-hexagon. In terms of Schur polynomials the result is
\begin{equation}\label{5.1}
\sum_{\mu\,;\,\ell(\mu)\le c} s_{\mu'}(1^a)s_{\mu'}(1^b)=
\prod_{i=1}^a\prod_{j=1}^b\prod_{k=1}^c \frac{i+j+k-1}{i+j+k-2},
\end{equation}
where the right hand side is MacMahon's formula.
(Here $1^a$ means $(1,\dots,1)$ with $a$ components.)
Comparing this
formula with the formula (\ref{4.6}) we find
\begin{equation}\label{5.2}
(-1)^{a+b} D_n(\prod_{\ell=1}^a(1-e^{-i\theta})
\prod_{\ell=1}^b(1-e^{i\theta}))=\prod_{i=1}^a\prod_{j=1}^b\prod_{k=1}^c
\frac{i+j+k-1}{i+j+k-2}.
\end{equation}
It has been conjectured by Keating and Snaith, \cite{KeSn}, that the
following result should hold for the moments of Riemann's
$\zeta$-function on the critical line,
\begin{equation}\label{5.3}
\lim_{T\to\infty}\frac 1{(\log T)^{k^2}}\frac 1{T} \int_0^T
|\zeta(1/2+it)|^{2k} dt=f_{\text{CUE}}(k)a(k),
\end{equation}
where $a(k)$ is a constant depending on the primes,
\begin{equation}\label{5.4}
f_{\text{CUE}}(k)=\lim_{n\to\infty} \frac 1{n^{k^2}}\int_{U(n)}
|Z(U,\theta)|^{2k} dU=\prod_{j=0}^{k-1}\frac {j!}{(j+k)!}
\end{equation}
and $Z(U,\theta)=\det(I-Ue^{-i\theta})$ is the characteristic
polynomial of the unitary matrix $U$. If we take $a=b=k$ in
(\ref{5.2}) and use (\ref{4.10}) we find
\begin{equation}\label{5.5}
\int_{U(n)} |Z(U,\theta)|^{2k} dU=
\prod_{i=1}^k\prod_{j=1}^k\prod_{\ell=1}^n
\frac{i+j+\ell-1}{i+j+\ell-2}=\prod_{j=0}^{n-1} \frac{j! (j+2k)!}{(j+k)!^2}
\end{equation}
as computed in \cite{KeSn} by different methods. Letting $n\to\infty$
we obtain the last expression in (\ref{5.4}). Hence, we see that the formula
(\ref{5.5}) has a curious combinatorial interpretation via MacMahon's
formula.

\label{lastpage}

\end{document}